\documentclass{amsart}

\def \Z{\mathbb Z}
\def \C{\mathbb C}
\def \D{{\mathcal{D}}}

\def \g{\mathfrak{g}}
\def \gl{\mathfrak{gl}}

\def \PAut{{\rm PAut}}
\def \Der{{\rm Der}}
\def \PDer{{\rm PDer}}

\def \Res{{\rm Res}}
\def \End{{\rm End}}
\def \PEnd{{\rm PEnd}}
\def \Hom{{\rm Hom}}
\def \<{\langle} 
\def \>{\rangle}
\def \be{\begin{equation}\label}
\def \ee{\end{equation}}
\def \bex{\begin{exa}\label}
\def \eex{\end{exa}}
\def \bl{\begin{lem}\label}
\def \el{\end{lem}}
\def \bt{\begin{thm}\label}
\def \et{\end{thm}}
\def \bp{\begin{prop}\label}
\def \ep{\end{prop}}
\def \br{\begin{rem}\label}
\def \er{\end{rem}}
\def \bc{\begin{coro}\label}
\def \ec{\end{coro}}
\def \bd{\begin{de}\label}
\def \ed{\end{de}}

\newtheorem{thm}{Theorem}[section]
\newtheorem{prop}[thm]{Proposition}
\newtheorem{coro}[thm]{Corollary}

\newtheorem{exa}[thm]{Example}
\newtheorem{lem}[thm]{Lemma}
\newtheorem{rem}[thm]{Remark}
\newtheorem{de}[thm]{Definition}

\makeatletter
\@addtoreset{equation}{section}

\makeatother
\makeatletter

\theoremstyle{definition}

\theoremstyle{remark}

\numberwithin{equation}{section}

\begin{document}

\title{Pseudoderivations, pseudoautomorphisms and 
simple current modules for vertex algebras}

\author{Haisheng Li}
\address{Department of Mathematics, Harbin Normal University, Harbin,
China}

\curraddr{Department of Mathematical Sciences, Rutgers University, 
Camden, NJ 08102}

\email{hli@camden.rutgers.edu}
\thanks{The author was supported in part by an NSA grant.}






\begin{abstract}
We exhibit a connection between 
Etingof-Kazhdan's notion of pseudoderivation and a certain construction of
simple current modules for a vertex operator algebra and meanwhile we
introduce and study a notion of pseudoendomorphism (pseudoautomorphism).
\end{abstract}

\maketitle


\section{Introduction}
Vertex (operator) algebras are analogous to classical Lie algebras and
associative algebras in many aspects, though there are also
significant differences in many aspects.  For a vertex (operator)
algebra $V$, one has the basic notions of derivation, endomorphism and
automorphism.  A derivation of a vertex algebra $V$ is a linear
endomorphism $d$ of $V$ such that
$$dY(u,x)v=Y(u,x)dv+Y(du,x)v\;\;\;\mbox{ for }u,v\in V.$$ Just as with
Lie algebras and associative algebras, 
the exponential $e^{d}$ of a locally finite
derivation $d$ of $V$ is an automorphism of $V$.  

In a study on formal deformation of vertex operator algebras, Etingof
and Kazhdan in \cite{ek} introduced and studied a notion of
pseudoderivation for a vertex algebra.  By definition, 
a pseudoderivation of a vertex algebra $V$ is a linear map $a(z)$
from $V$ to $\C((z))\otimes V$, where $z$ is a formal variable, such
that
\begin{eqnarray}
& &[L(-1),a(z)]=-\frac{d}{dz}a(z),\label{eintro-pseuoder1}\\
& &[a(z),Y(v,x)]=Y(a(z+x)v,x)\;\;\;\;\mbox{ for }v\in V.\label{eintro-pseuoder2}
\end{eqnarray}
It was proved therein that $a(z){\bf 1}=0$.
Just as derivations are important in the study of formal deformations
of classical Lie algebras, pseudoderivations are important in the
study of formal deformations of vertex operator algebras.

In another (seemingly unrelated) study \cite{li-super},
motivated by the physics superselection principle in 
(algebraic) conformal field theory,
we studied how to deform the adjoint module 
for a vertex operator algebra $V$ to obtain new modules.
The formulation of the physics superselection principle 
is based on the following simple fact: For an algebra of many types
(Lie algebras, associative algebras and vertex algebras),
if $\sigma$ is an endomorphism of the algebra, then for any
module $U$, we have a ``new'' module structure on $U$ defined by
$a*u=\sigma(a)u$ for $a\in R,\; u\in U$. 
Having noticed that the $\sigma$-twisting of
the adjoint representation by an automorphism $\sigma$
is always isomorphic to the adjoint representation we considered
twisting by deformed endomorphisms of a certain type.
It was proved in \cite{li-super} that if
$\Delta(z)\in \Hom (V,\C((z))\otimes V)$ satisfies the conditions
\begin{eqnarray}
& &\Delta(z){\bf 1}={\bf 1}\;(=1\otimes {\bf 1}),\label{eintro-pseudoend1}\\
& &[L(-1),\Delta(z)]=-\frac{d}{dz}\Delta(z),\label{eintro-pseudoend2}\\
& & \Delta(z)Y(v,x)=Y(\Delta(z+x)v,x)\Delta(z)
\;\;\;\mbox{ for }v\in V,\label{eintro-pseudoend3}
\end{eqnarray}
then for any $V$-module $(W,Y_{W})$, $(W,Y_{W}(\Delta(x)\cdot,x))$ 
carries the structure of a $V$-module. 
(Notice that if $\Delta(z)=\sigma$ is independent of $z$, then
the above three conditions exactly amount to that 
$\sigma$ is an endomorphism of vertex algebra $V$.)
Furthermore, we gave an explicit construction
of such invertible $\Delta(z)$ and we applied the results
to the vertex operator algebras associated with affine Lie algebras,
with Heisenberg Lie algebras and with nondegenerate even lattices.
Indeed, the twisting of the adjoint module gives rise to
{\em new} modules, which were proved to be ``simple current modules''
in a certain sense. In this way we obtained a construction of
simple current modules.

Now, it is manifest that Etingof-Kazhdan's notion of pseudoderivation and 
our construction of simple current modules are much more than simply related.
The main purpose of this short note is to exhibit this connection explicitly
and conceptually. In Section 2, we formulate notions of
pseudoendomorphism and pseudoautomorphism, 
using (\ref{eintro-pseudoend1})-(\ref{eintro-pseudoend3}).
It is proved that that if $a(z)$ is a locally nilpotent pseudoderivation 
of $V$ such that $[a(z_{1}),a(z_{2})]=0$,
$e^{a(z)}$ is a pseudoautomorphism of $V$.
We show that the Lie algebra $\PDer (V)$ is canonically isomorphic to
the Lie algebra of all $\C((z))$-linear derivations of the tensor product vertex algebra
$\C((z))\otimes V$ (see \cite{flm}), where $\C((z))$ is 
the Borcherds' vertex algebra associated with the unital commutative
associative algebra $\C((z))$ together with the derivation 
$D=d/dz$. Similarly, we show that any endomorphism of
$\C((z))\otimes V$, restricted to $V$, is a pseudoendomorphism 
of $V$, and that the $\C((z))$-linear extension of
every pseudoendomorphism of $V$ is an endomorphism
of vertex algebra $\C((z))\otimes V$.
In Section 3, we reinterpret the construction of simple current modules
in terms of pseudoautomorphisms.

\section{Pseudoderivations and pseudoendomorphisms}
In this section we review Etingof-Kazhdan's notion of pseudoderivation
and define a notion of pseudoendomorphism
for a vertex algebra $V$ and we then relate pseudoderivations 
and pseudoendomorphism of $V$ with
derivations and endomorphisms of the tensor product vertex algebra 
$\C((z))\otimes V$.

We shall not recall the full definition of a vertex algebra,
but for convenience we just mention the main ingredients.
For a vertex algebra $V$, we have the vacuum vector ${\bf 1}$ and
the vertex operator map $Y:V\rightarrow \Hom (V,V((x)))\subset (\End V)[[x,x^{-1}]]$, where
$$Y(v,x)=\sum_{n\in \Z}v_{n}x^{-n-1}\;\;\;(\mbox{with }v_{n}\in \End V).$$
There is a canonical linear operator $\D$ on $V$ defined by
\begin{eqnarray}
\D (v)=v_{-2}{\bf 1}\;\;\;\mbox{ for }v\in V.
\end{eqnarray}
The basic properties about $\D$ are
\begin{eqnarray}
& &Y(u,x){\bf 1}=e^{x{\D}}u,\\
& &Y(u,x)v=e^{x{\D}}Y(v,-x)u,\\
& &[{\D},Y(v,x)]=Y({\D}v,x)=\frac{d}{dx}Y(v,x)
\;\;\;\mbox{ for }u,v\in V.
\end{eqnarray}

An {\em endomorphism} of a vertex algebra $V$ is a linear 
endomorphism $f$ of $V$ such that $f({\bf 1})={\bf 1}$ and 
\begin{eqnarray}
fY(u,x)v=Y(f(u),x)f(v)\;\;\;\mbox{ for }u,v\in V.
\end{eqnarray}
An automorphism of a vertex algebra $V$ is a bijective endomorphism.

A {\em derivation} of a vertex algebra $V$ (cf. [B], [Lia])
is a linear endomorphism $d$ such that
\begin{eqnarray}\label{ederivation-definition}
dY(u,x)v=Y(u,x)dv+Y(du,x)v\;\;\;\mbox{ for }u,v\in V.
\end{eqnarray}
The operator $\D$ is a derivation of $V$ and
it follows from Borcherds' commutator formula
that for any $v\in V$, $v_{0}$ is a derivation of $V$. 

If $d$ is a derivation of $V$, one can show
that $d{\bf 1}=0$ and $[{\D},d]=0$.  It is straightforward to see
that all the derivations of $V$ form a Lie subalgebra $\Der V$ of
$\gl(V)$ (the general Lie algebra), with an ideal consisting of
$v_{0}$ for $v\in V$. 

\br{rderivations} 
{\em If $d$ is a derivation of $V$, which is locally
finite in the sense that for any $v\in V$, the vectors $d^{n}w$ for $n\ge 0$
span a finite-dimensional subspace of $V$, then $e^{d}$ is an automorphism
of the vertex algebra $V$.} 
\er

The following notion is due to Etingof and Kazhdan \cite{ek}:

\bd{dpseudoderivation-ek}
{\em Let $V$ be a vertex algebra. 
A {\em pseudoderivation} of a vertex algebra $V$ is 
a linear homomorphism $a(z)$ from $V$ to $\C((z))\otimes V$ such that
\begin{eqnarray}
& &[{\D},a(z)]=-\frac{d}{dz}a(z),\label{edefinition-pseu1}\\
& &[a(z),Y(u,x)]=Y(a(z+x)u,x)\;\;\;\mbox{ for }u\in V,
\label{edefinition-pseu2}
\end{eqnarray}
where $\D$ acts on $\C((z))\otimes V$ as $1\otimes \D$ and
the map $Y$ is $\C((z))[[x]]$-linearly extended to $\C((z))[[x]]\otimes V$.}
\ed

All the pseudoderivations of vertex algebra $V$ clearly form a
$\C$-subspace of $\Hom (V,\C((z))\otimes V)$, which is denoted by 
$\PDer (V)$, as in \cite{ek}.   Note that we have
\begin{eqnarray}
\Hom _{\C}(V,\C((z))\otimes V)\simeq \End_{\C((z))}(\C((z))\otimes V)
\end{eqnarray}
as $\C$-vector spaces.
As it was pointed out in \cite{ek}, $\PDer (V)$ is naturally
a Lie subalgebra of $\gl(\C((z))\otimes V)$.

We have the following results, among which
the second assertion was proved in \cite{ek}:

\bp{p-pseudoderivation-red}
In Definition \ref{dpseudoderivation-ek}, (\ref{edefinition-pseu1})
follows from (\ref{edefinition-pseu2}). Furthermore, for any pseudoderivation $a(z)$
we have $a(z){\bf 1}=0$.
\ep

\begin{proof}
Since $Y({\bf 1},x)=1$, from (\ref{edefinition-pseu2}) we have
$Y(a(z+x){\bf 1},x)=0$.
As the map $Y$ is injective (see [FHL]), 
we have $a(z+x){\bf 1}=0$. Thus $a(z){\bf 1}=0$. 

Now applying both sides of (\ref{edefinition-pseu2}) to the vacuum
vector ${\bf 1}$ and using the fact that $a(z){\bf 1}=0$, we get
\begin{eqnarray*}
a(z)Y(u,x){\bf 1}=Y(a(z+x)u,x){\bf 1}.
\end{eqnarray*}
Thus
\begin{eqnarray*}
a(z)e^{x\D}u=e^{x\D}a(z+x)u=e^{x\D}e^{xd/dz}a(z)u.
\end{eqnarray*}
That is,
\begin{eqnarray}
e^{-x\D}a(z) e^{x\D}u=e^{xd/dz}a(z)u,
\end{eqnarray}
{}from which we immediately have (\ref{edefinition-pseu1}).
\end{proof}

The following result, due to \cite{ek}, gives a construction of
pseudoderivations:

\bp{pek}
Let $V$ be a vertex algebra, let $v\in V$ and let $f(z)\in \C((z))$. Then
\begin{eqnarray}
X_{f,v}=\Res_{x} f(z+x)Y(v,x)=\sum_{n\ge 0}\frac{f^{(n)}(z)}{n!}v_{n}
\end{eqnarray}
is a pseudoderivation of $V$.
\ep

As we have done so in Definition \ref{dpseudoderivation-ek}, we extend 
the $\C$-linear map $Y$ from $V$ to $\Hom (V,V((x)))$
to a $\C((z))$-linear map from $\C((z))\otimes V$ to
$(\End_{\C((z))} (\C((z))\otimes V))[[x,x^{-1}]]$.
With $\C((z))$ as a field, $\C((z))\otimes V$ equipped with 
the $\C((z))$-linear map $Y$ is naturally a vertex algebra over $\C((z))$.

On the other hand,  $\C((z))$ is a unital commutative associative algebra over $\C$
with $d/dz$ as a derivation. {}From [B],
$\C((z))$ is a vertex algebra with $1$ as the vacuum vector where
\begin{eqnarray}
Y(f(z),x)g(z)=f(z+x)g(z)=(e^{xd/dz}f(z))g(z)
\end{eqnarray}
for $f(z),g(z)\in \C((z))$. Then for any vertex algebra $V$, from [FHL], 
$\C((z))\otimes V$ is a vertex algebra where the vertex operator map
$\widehat{Y}$ is given by
\begin{eqnarray}
\widehat{Y}(f(z)\otimes u,x)(g(z)\otimes v)
=f(z+x)g(z)\otimes Y(u,x)v
\end{eqnarray}
for $f(z),g(z)\in \C((z)),\; u,v\in V$. 
Throughout this note, we always refer the vertex algebra
$\C((z\C((z))\otimes V$ to this specific vertex algebra structure.

In the following we shall use both maps $Y$ and 
$\widehat{Y}$ on $\C((z))\otimes V$. 
It is very important for us to distinguish them.
For any vector $w(z)\in \C((z))\otimes V$, we have
\begin{eqnarray}\label{etensor-vertex-operator-map}
\widehat{Y}(w,x)=Y(w(z+x),x)
=e^{x\frac{\partial}{\partial z}}Y(w(z),x).
\end{eqnarray}

Denote by $\Der_{\C((z))}(\C((z))\otimes V)$ the space of $\C((z))$-linear derivations
of vertex algebra $\C((z))\otimes V$.

\bp{ppseudoderivation-charc}
Let $V$ be any vertex algebra (over $\C$). Every derivation of
vertex algebra $\C((z))\otimes V$, restricted to $V$, is a pseudoderivation of $V$.
On the other hand, the $\C((z))$-linear extension of any pseudoderivation of $V$
is a derivation of $\C((z))\otimes V$.
Furthermore, we have
\begin{eqnarray}\label{eidentification}
\Der_{\C((z))} (\C((z))\otimes V)=\PDer(V).
\end{eqnarray}
\ep

\begin{proof}
Let $\sigma$ be a derivation of the vertex algebra $\C((z))\otimes V$.
By restricting (the domain of) $\sigma$ to $V$ we have
a $\C$-linear map $a(z)$ from $V$ to $\C((z))\otimes V$. 
Using (\ref{etensor-vertex-operator-map}) we get
\begin{eqnarray*}
a(z) Y(u,x)v-Y(u,x)a(z)v
&=&a(z) \widehat{Y}(u,x)v-\widehat{Y}(u,x)a(z)v\\
&=&\widehat{Y}(a(z)u,x)v\\
&=&Y(a(z+x)u,x)v.
\end{eqnarray*}
In view of Proposition \ref{p-pseudoderivation-red}, 
$a(z)$ is a pseudoderivation of $V$.

On the other hand, let $a(z)$ be a pseudoderivation of $V$.
Consider $a(z)$ as a $\C((z))$-linear endomorphism of $\C((z))\otimes V$.
For $f(z),g(z)\in \C((z)),\; u,v\in V$, using (\ref{etensor-vertex-operator-map}) we have
\begin{eqnarray*}
& &a(z)\widehat{Y}(f(z)\otimes u,x)(g(z)\otimes v)\nonumber\\
&=&a(z)(f(z+x)g(z)\otimes Y(u,x)v)\nonumber\\
&=&f(z+x)g(z) a(z)Y(u,x)v\nonumber\\
&=&f(z+x)g(z)Y(u,x)a(z)v+f(z+x)g(z)Y(a(z+x)u,x)v\nonumber\\
&=&Y(f(z+x)u,x)g(z)a(z)v+Y(f(z+x)a(z+x)u,x)g(z)v\nonumber\\
&=&\widehat{Y}(f(z)u,x)a(z)(g(z)v)+
\widehat{Y}(a(z)f(z)u,x)g(z)v.
\end{eqnarray*}
This proves that $a(z)$ is a derivation of $\C((z))\otimes V$.
The identification of the two spaces in (\ref{eidentification}) is clear.
\end{proof}

\br{rdiff-proof1}
{\em For any $f(z)\in \C((z)),\; v\in V$, $(f(z)v)_{0}$ is a derivation
of vertex algebra $\C((z))\otimes V$. We have
\begin{eqnarray}
(f(z)v)_{0}=\Res_{x}\widehat{Y}(f(z)v,x)=\Res_{x}f(z+x)Y(v,x)
=X_{f,v}.
\end{eqnarray}
In view of Proposition \ref{ppseudoderivation-charc},
$X_{f,v}$ is a pseudoderivation of $V$.
This gives a (slightly) different proof of Proposition \ref{pek}.}
\er

Define a $\C$-linear map
\begin{eqnarray}
X:\  \C((z))\otimes V\rightarrow \PDer (V),\ \  (f,v)\mapsto X_{f,v}.
\end{eqnarray}
Denote by $\hat{\D}$ the $\D$-operator of $\C((z))\otimes V$.
We have
\begin{eqnarray}
\hat{\D}=d/dz\otimes 1+1\otimes \D.
\end{eqnarray}
Since $(\hat{\D} (f\otimes v))_{0}=0$ on $\C((z))\otimes V$ 
for any $f\in \C((z)),\; v\in V$, we have
$\hat{\D} (\C((z))\otimes V)\subset \ker X$.
Then we have a linear map, also denoted by $X$, from 
$(\C((z))\otimes V)/\hat{\D}(\C((z))\otimes V)$ to $\PDer (V)$.
Recall from \cite{bor} that for any vertex algebra $U$,  $U/\D U$ is a Lie algebra with
\begin{eqnarray}
[a+\D U,b+\D U]=a_{0}b+\D U\;\;\;\mbox{ for }a,b\in U.
\end{eqnarray}
Now, taking $U=\C((z))\otimes V$, we have a Lie algebra
$(\C((z))\otimes V)/\hat{\D}(\C((z))\otimes V)$ with
\begin{eqnarray}
[f(z)\otimes u, g(z)\otimes v]=
(f(z)\otimes u)_{0}(g(z)\otimes v)
=\sum_{n\ge 0}\frac{1}{n!}f^{(n)}(z)g(z)\otimes u_{n}v
\end{eqnarray}
for $f(z),g(z)\in \C((z)),\; u,v\in V$.
On the other hand, we have
\begin{eqnarray*}
[X_{f,u},X_{g,v}]&=&\sum_{m,n\ge 0}\frac{1}{m!}\frac{1}{n!}f^{(m)}g^{(n)}[u_{m},v_{n}]\\
&=&\sum_{m,n\ge 0}\sum_{i\ge 0}\binom{m}{i}\frac{1}{m!}\frac{1}{n!}f^{(m)}g^{(n)}(u_{i}v)_{m+n-i}\\
&=&\sum_{m',n\ge 0}\sum_{i\ge 0}\binom{m'+i}{i}\frac{1}{(m'+i)!}\frac{1}{n!}f^{(m'+i)}g^{(n)}(u_{i}v)_{m'+n}\\
&=&\sum_{k\ge 0}\sum_{n=0}^{k}\sum_{i\ge 0}\binom{k-n+i}{i}
\frac{1}{(k-n+i)!}\frac{1}{n!}f^{(k-n+i)}g^{(n)}(u_{i}v)_{k}\\
&=&\sum_{k\ge 0}\sum_{i\ge 0}
\frac{1}{k!}\frac{1}{i!}(f^{(i)}g)^{(k)}(u_{i}v)_{k}\\
&=&X_{[f\otimes u,g\otimes v]}.
\end{eqnarray*}
Thus $X$ is a Lie algebra homomorphism.

To summarize we have:

\bp{plie-hom}
For any vertex algebra $V$, the linear map 
\begin{eqnarray}
X: & &(\C((z))\otimes V)/\hat{\D}(\C((z))\otimes V)\rightarrow \PDer (V),\nonumber\\
& &(f(z),v)\mapsto X_{f,v}
\end{eqnarray}
is a Lie algebra homomorphism.
\ep

Next, we study notions of pseudoendomorphism
and pseudoautomorphism.

\bd{dpseudoendomorphism}
{\em A {\em pseudoendomorphism} of a vertex algebra $V$ is a 
$\C$-linear map $A(z)$ from $V$ to $\C((z))\otimes V$ such that}
\begin{eqnarray}
& &A(z) {\bf 1}={\bf 1}\;(=1\otimes {\bf 1}),\label{epseudoend1}\\
& &[{\D},A(z)]=-\frac{d}{dz}A(z),\label{epseudoend2}\\
& &A(z)Y(v,x)=Y(A(z+x)v),x)A(z)\;\;\;\mbox{ for }v\in V.\label{epseudo-end-3}
\end{eqnarray}
\ed

A {\em pseudoautomorphism} of a vertex algebra $V$ is an invertible
pseudoendomorphism of $V$.  Clearly, all the pseudoautomorphisms of $V$
form a group, which we denote by $\PAut (V)$.

\br{rpseudoendomorphism}
{\em It is clear that a constant pseudoendomorphism of $V$
exactly amounts to an endomorphism of vertex algebra $V$.}
\er

The following result, analogous to Proposition \ref{p-pseudoderivation-red},
shows that condition (\ref{epseudoend2})
in Definition \ref{dpseudoendomorphism} is redundant:

\bl{lredundant} 
Let $V$ be a vertex algebra and let $A(z)\in \Hom
(V,\C((z))\otimes V)$ be such that (\ref{epseudoend1}) and
(\ref{epseudo-end-3}) hold.  Then $A(z)$ is a pseudoendomorphism of
$V$.  
\el

\begin{proof} 
Applying (\ref{epseudo-end-3}) to ${\bf 1}$ and 
using the assumption that $A(z){\bf 1}={\bf 1}$ we get
\begin{eqnarray}
A(z)Y(v,x){\bf 1}=Y(A(z+x)v,x){\bf 1}
=e^{x\frac{\partial}{\partial z}}Y(A(z)v,x){\bf 1}.
\end{eqnarray}
Just as in the second paragraph of the proof of 
Proposition \ref{p-pseudoderivation-red}, we get
$$[A(z),{\D}]v=\frac{d}{dz}A(z)v,$$
proving (\ref{edelta-property2}). 
Thus $A(z)$ is a pseudoendomorphism of $V$.
\end{proof}

We also have the following analogue of Proposition \ref{ppseudoderivation-charc}:

\bp{ppseudoautomorphism}
Let $V$ be a vertex algebra. Every endomorphism of vertex algebra
$\C((z))\otimes V$, restricted to $V$, is a pseudoendomorphism of 
vertex algebra $V$. On the other hand,
the $\C((z))$-linear extension of any pseudoendomorphism of 
vertex algebra $V$ is an endomorphism of vertex algebra
$\C((z))\otimes V$. Furthermore, we have
\begin{eqnarray}
\PEnd (V)=\End_{\C((z))} (\C((z))\otimes V),
\end{eqnarray}
where $\End_{\C((z))} (\C((z))\otimes V)$ stands for the space
of all the $\C((z))$-linear endomorphisms of vertex algebra
$\C((z))\otimes V$.
\ep

\begin{proof} Let $\sigma$ be an endomorphism of vertex algebra
$\C((z))\otimes V$. Denote by $A_{\sigma}(z)$ the restriction
of $\sigma$ as a $\C$-linear map from $V$ to $\C((z))\otimes V$.
Then 
\begin{eqnarray*}
A_{\sigma}(z){\bf 1}=\sigma(1\otimes {\bf 1})
=1\otimes {\bf 1}={\bf 1}
\end{eqnarray*}
and for $u,v\in V$, using (\ref{etensor-vertex-operator-map}) we get
\begin{eqnarray*}
A_{\sigma}(z)Y(u,x)v=
\sigma \widehat{Y}(u,x)v
&=&\widehat{Y}(\sigma(u),x)\sigma(v)\nonumber\\
&=&\widehat{Y}(A_{\sigma}(z)u,x)A_{\sigma}(z)v\nonumber\\
&=&Y(A_{\sigma}(z+x)u,x)A_{\sigma}(z)v.
\end{eqnarray*}
In view of Lemma \ref{lredundant}, $A_{\sigma}(z)$ is a pseudoendomorphism of 
vertex algebra $V$.

On the other hand, let $A(z)$ be a pseudoendomorphism of 
vertex algebra $V$. Denote by $\bar{A}(z)$ the corresponding $\C((z))$-linear
endomorphism of $\C((z))\otimes V$.
Then for $f(z)\in \C((z)),\;v\in V$, we have
\begin{eqnarray}
& &\bar{A}(z) \widehat{Y}( f(z)u,x)(g(z)v)\nonumber\\
&=&\bar{A}(z) \left(f(z+x)g(z)Y(u,x)v\right)\nonumber\\
&=&f(z+x)g(z)A(z)Y(u,x)v\nonumber\\
&=&f(z+x)g(z)Y(A(z+x)u,x)A(z)v\nonumber\\
&=&Y(f(z+x)A(z+x)u,x)g(z) A(z)v\nonumber\\
&=&Y(\bar{A}(z+x)f(z+x)u,x)g(z) A(z)v\nonumber\\
&=&\widehat{Y}(\bar{A}(z) (f(z)u),x)\bar{A}(z)(g(z)v).
\end{eqnarray}
That is, $\bar{A}(z)$ is an endomorphism of vertex algebra $\C((z))\otimes V$. 
\end{proof}

Just as the exponential of a locally nilpotent derivation of a Lie algebra is
an automorphism we have:

\bp{pexponentiation}
Let $V$ be a vertex algebra and let
$a(z)$ be a pseudoderivation of $V$ such that $[a(z),a(x)]=0$ and such that 
for every $v\in V$, there exists a positive integer $k$ such that
$a(z)^{k}v=0$. Then $\exp a(z)$ is a pseudoautomorphism of $V$.
\ep

\begin{proof}
It is clear that under the assumptions, $\exp a(z)$ is 
a well defined linear map from $V$ to $\C((z))\otimes V$. 
As $a(z){\bf 1}=0$, we have $(\exp a(z)){\bf 1}={\bf 1}$.
Since
$[a(z),Y(v,x)]=Y(a(z+x)v,x)$ for $v\in V$,
it follows immediately that
\begin{eqnarray}
(\exp a(z)) Y(v,x)\exp (-a(z))=Y( (\exp a(z+x))v,x).
\end{eqnarray}
By Lemma \ref{lredundant},  $\exp a(z)$ is a pseudoautomorphism of $V$.
\end{proof}

\bex{exalattice}
{\em Let $V=L_{\hat{\g}}(\ell,0)$ be the simple vertex operator algebra
associated with an affine Kac-Moody Lie algebra $\hat{\g}$
and with a positive integer $\ell$. Then $L_{\hat{\g}}(\ell,0)$ is an integrable
$\hat{\g}$-module (see \cite{kac}). Let $u$ be a root vector of $\g$. We have
$[u(x_{1}),u(x_{2})]=0$ and for every $n\in \Z$,
$u(n)$ acts locally nilpotently on  $L_{\hat{\g}}(\ell,0)$.
Furthermore, for any $w\in L_{\hat{\g}}(\ell,0)$ we have
$u(n)w=0$ for $n$ sufficiently large.
It follows that for any $f(z)\in \C((z))$ and for any $w\in  L_{\hat{\g}}(\ell,0)$, 
there exists a nonnegative integer $k$ such that
$X_{f,u}^{k}w=0$.
(Recall that $X_{f,u}=\sum_{n\ge 0}f^{(n)}(z)u(n)$.)
By Proposition \ref{pexponentiation}, we have a pseudoautomorphism $\exp X_{f,u}$.
It will be interesting to study the subgroup of $\PAut (V)$,
generated by all such $\exp X_{f,u}$.}
\eex

\section{Simple current modules and pseudoautomorphisms}
In this section we reveal a connection between
pseudoautomorphisms and simple current modules for a vertex algebra.

Let us first recall the following results from \cite{li-super}:

\bp{pmoduledeformation}
Let $V$ be a vertex algebra and let 
$\Delta(z)\in \Hom (V,\C((z))\otimes V)$. If 
$(V,Y(\Delta(x)\cdot,x))$ carries the structure of a $V$-module, 
then
\begin{eqnarray}
& &\Delta(z){\bf 1}={\bf 1}\;(1\otimes {\bf 1}),\label{edelta-property1}\\
& &[{\D},\Delta(z)]=-\frac{d}{dz}\Delta(z),\label{edelta-property2}\\
& &\Delta(z)Y(v,x)=Y(\Delta(z+x)v,x)\Delta(z)\;\;\;\mbox{for }v\in V. 
\label{edelta-property3}
\end{eqnarray}
On the other hand, if 
$\Delta(z)\in \Hom (V,\C((z))\otimes V)$ satisfies 
the above three conditions, then for any $V$-module $(W,Y_{W})$,
$(W,Y_{W}(\Delta(x)\cdot,x))$ carries the structure of a $V$-module.
\ep

The following notion is due to \cite{sy}:

\bd{dsimplecurrent}
{\em Let $V$ be a simple vertex operator algebra. 
An irreducible (simple) $V$-module $W$ is called 
a {\em simple current $V$-module} 
if the matrix of the left multiplication associated to
the equivalence class of $W$ in the fusion algebra of $V$
is a permutation.}
\ed

\bp{psimplecurrent}
Let $V$ be an irreducible vertex algebra and let 
$\Delta(z)\in \Hom (V,\C((z))\otimes V)$ such that the three conditions 
in Proposition \ref{pmoduledeformation} hold. Suppose that $\Delta(z)$ 
is invertible in the sense that there exists $A(z)\in \Hom
(V,\C((z))\otimes V)$ such that $A(z)\Delta(z)=\Delta(z)A(z)=1$. 
Then $(V,Y(\Delta(x)\cdot,x))$ carries 
the structure of a simple current $V$-module.
\ep

\br{rabout-deltaz}
{\em In \cite{li-super}, it was assumed that $\Delta(z)\in \Hom
(V,V[z,z^{-1}])$, however, all the proofs just work fine if  
$\Delta(z)\in \Hom (V,\C((z))\otimes V)$.}
\er

The following result was obtained in \cite{li-super} (see also \cite{li-twisted}):

\bp{pli}
Let $V$ be a vertex algebra and let $h\in V$ be such that
\begin{eqnarray}
h_{1}h\in \C {\bf 1}\;\;\mbox{ and }\;\;h_{n}h=0\;\;\;\mbox{ for }n\ge
0,\; n\ne 1
\end{eqnarray}
and such that $h_{0}$ acts semisimply on $V$ with integral eigenvalues
and $h_{n}$ for $n\ge 1$ act locally nilpotently on $V$.
Set
\begin{eqnarray}
\Delta(h,z)=z^{h(0)}\exp\left(\sum_{n\ge 1}\frac{h(n)}{-n}(-z)^{-n}\right).
\end{eqnarray}
Then $\Delta(h,z)$ satisfies the conditions 
(\ref{edelta-property1})-(\ref{edelta-property3}).
\ep

Now, in terms of our notions of pseudoendomorphism and pseudoautomorphism
Propositions \ref{pmoduledeformation} and \ref{psimplecurrent}
can be reformulated as follow:

\bp{pmodulede-pseudo-endo}
Let $V$ be a vertex algebra and let 
$\Delta(z)\in \Hom (V,\C((z))\otimes V)$. If 
$(V,Y(\Delta(x)\cdot,x))$ carries the structure of a $V$-module, 
then $\Delta(z)$ is a pseudoendomorphism of $V$.
On the other hand, if 
$\Delta(z)$ is a pseudoendomorphism of $V$, then for any $V$-module $(W,Y_{W})$,
$(W,Y_{W}(\Delta(x)\cdot,x))$ carries the structure of a $V$-module.
Furthermore, if $\Delta(z)$ is invertible, $(V,Y(\Delta(x)\cdot,x))$ is a
simple current $V$-module.
\ep
 
Next, we discuss certain connections between $V$-modules and
$\C((z))\otimes V$-modules. Since $V$ is a vertex subalgebra of
$\C((z))\otimes V$, any $\C((z))\otimes V$-module is naturally
a $V$-module. On the other hand, any $V$-module can be naturally
made a $\C((z))\otimes V$-module by the following result:

\bp{paffine-module}
Let $(W,Y_{W})$ be a $V$-module. Define a $\C$-linear map
$$\widehat{Y}_{W}: \C((z))\otimes V\rightarrow (\End W)[[x,x^{-1}]]$$
by
\begin{eqnarray}
\widehat{Y}_{W}(f(z)\otimes v,x)=f(x)Y_{W}(v,x)
\end{eqnarray}
for $f(z)\in \C((z)),\; v\in V$.
Then $(W,\widehat{Y}_{W})$ carries the structure of a $\C((z))\otimes V$-module.
\ep

\begin{proof} It is clear that $\widehat{Y}_{W}(f(z)\otimes v,x)\in \Hom (W,W((x)))$
for $f(z)\in \C((z)),\; v\in V$.
We have
$$\widehat{Y}_{W}(1\otimes {\bf 1},x)=Y_{W}({\bf 1},x)=1_{W}.$$
Furthermore, for $f(z),g(z)\in \C((z)),\; u,v\in V$, we have
\begin{eqnarray*}
& &\widehat{Y}_{W}(f(z)\otimes u,x_{1})=f(x_{1})Y_{W}(u,x_{1}),\\
& &\widehat{Y}_{W}(g(z)\otimes v,x_{2})=g(x_{2})Y_{W}(v,x_{2}),
\end{eqnarray*}
and 
\begin{eqnarray*}
& &x_{2}^{-1}\delta\left(\frac{x_{1}-x_{0}}{x_{2}}\right) 
\widehat{Y}_{W}(\widehat{Y}(f(z)\otimes u,x_{0})(g(z)\otimes v,x_{2})\\
&=&x_{2}^{-1}\delta\left(\frac{x_{1}-x_{0}}{x_{2}}\right) 
\widehat{Y}_{W}(f(z+x_{0})g(z)\otimes Y(u,x_{0})v,x_{2})\\
&=&x_{2}^{-1}\delta\left(\frac{x_{1}-x_{0}}{x_{2}}\right) 
f(x_{2}+x_{0})g(x_{2})Y_{W}(Y(u,x_{0})v,x_{2})\\
&=&x_{2}^{-1}\delta\left(\frac{x_{1}-x_{0}}{x_{2}}\right) 
f(x_{1})g(x_{2})Y_{W}(Y(u,x_{0})v,x_{2}).
\end{eqnarray*}
It follows immediately that the Jacobi identity holds.
Thus $(W,\widehat{Y}_{W})$ carries the structure 
of a $\C((z))\otimes V$-module.
\end{proof}

\br{rnewproof}
{\em We here give a new proof for the assertion in Proposition \ref{pmodulede-pseudo-endo}
that if $\Delta(z)$ is a pseudoendomorphism of $V$, for any
$V$-module $(W,Y_{W})$, $(W,Y_{W}(\Delta(x)\cdot,x))$
is a $V$-module.
In view of Proposition \ref{ppseudoautomorphism}, a pseudoendomorphism $\Delta(z)$ of $V$ 
is naturally an endomorphism of 
the tensor product vertex algebra $\C((z))\otimes V$.
For any $V$-module $(W,Y_{W})$, by Proposition \ref{paffine-module},
we have a $\C((z))\otimes V$-module
$(W,\widehat{Y}_{W})$. By the simple fact,
$(W,\widehat{Y}_{W}\circ \Delta(z))$ is also a 
$\C((z))\otimes V$-module. In particular,
$(W,\widehat{Y}_{W}\circ \Delta(z))$ is a $V$-module.
Furthermore, for $u\in V$, we have
\begin{eqnarray}
\widehat{Y}_{W}(\Delta(z)u,x)
=Y_{W}(\Delta(x)u,x).
\end{eqnarray}
This shows that $(W,Y_{W}(\Delta(x)\cdot,x))$ is a $V$-module.}
\er

\br{rfuture}
{\em Recall from Proposition \ref{pli} the pseudoautomorphism $\Delta(h,z)$.
Since $[h(m),h(n)]=0$ for $m,n\ge 0$, we can write $\Delta(h,z)$ as
\begin{eqnarray}
\Delta(h,z)=\exp\left(h(0)\log z+\sum_{n\ge 1}(-1)^{n-1}\frac{h(n)}{n}z^{-n}\right).
\end{eqnarray}
Notice that if we could take $f(z)=\log z$ in the formula
$$X_{f,v}=\sum_{n\ge 0}\frac{1}{n!}f^{(n)}(z)v_{n}=\Res_{x}e^{x(d/dz)}f(z) Y(v,x),$$
we would have
$$X_{\log z,v}= \sum_{n\ge 0} \frac{1}{n!}(d/dz)^{n}(\log z) v_{n}
=v_{0}\log z +\sum_{n\ge 1} (-1)^{n-1}\frac{1}{n}v_{n}z^{-n}.$$
Now, $X_{\log z,v}$ is a $\C$-linear map from $V$ to $\C((z))[\log z]\otimes V$.
One can easily proves
\begin{eqnarray}
[X_{\log z,v},Y(u,x)]=Y(e^{x(d/dz)}X_{\log z,v}u,x)\;\;\;\mbox{ for }u\in V.
\end{eqnarray}
This suggests us to study ``logarithmic pseudoderivations'' of $V$, which are
$\C$-linear maps $a(z)$ from $V$ to $\C((z))[\log z]\otimes V$
satisfying the condition
\begin{eqnarray}
[a(z),Y(u,x)]=Y(e^{x(d/dz)}a(z)u,x)\;\;\;\mbox{ for }u\in V.
\end{eqnarray}
One can also define
a notion of logarithmic pseudoendomorphism of $V$ as a linear map
$A(z)$ from $V$ to $\C((z))[\log z]\otimes V$ such that
\begin{eqnarray}
A(z)Y(u,x)v=Y(e^{x(d/dz)}A(z)u,x)A(z)v\;\;\;\mbox{ for }u,v\in V.
\end{eqnarray}
We shall study this in a future publication.}
\er

\end{document}